\documentclass[aos,nameyear,dvips]{arximspdf}

\doi{10.1214/10-AOS801}
\volume{38}
\issue{5}
\pubyear{2010}
\firstpage{3191}
\lastpage{3216}

\makeatletter

\newcommand{\interleave}{\Vert\hspace*{-1.4pt}\vert}

\newtheorem{theorem}{Theorem}[section]
\newtheorem{corollaryNumThm}[theorem]{Corollary}

\makeatother

\begin{document}
\begin{frontmatter}

\title{On information plus noise kernel random matrices}
\runtitle{Information plus noise kernel random matrices}

\begin{aug}
\author[A]{\fnms{Noureddine} \snm{El Karoui}\corref{}\thanksref{t1}\ead[label=e1]{nkaroui@stat.berkeley.edu}}
\runauthor{N. El Karoui}
\affiliation{University of California, Berkeley}
\address[A]{Department of Statistics\\
University of California, Berkeley\\
367 Evans Hall\\
Berkeley, California 94720-3860\\
USA\\
\printead{e1}} 
\end{aug}

\thankstext{t1}{Supported by the NSF Grants DMS-06-05169,
DMS-08-47647 (CAREER) and an Alfred P. Sloan research fellowship.}

\received{\smonth{7} \syear{2009}}
\revised{\smonth{12} \syear{2009}}

%
\begin{abstract}
Kernel random matrices have attracted a lot of interest in recent
years, from both practical and theoretical standpoints. Most of the
theoretical work so far has focused on the case were the data is
sampled from a low-dimensional structure. Very recently, the first
results concerning kernel random matrices with high-dimensional input
data were obtained, in a setting where the data was sampled from a
genuinely high-dimensional structure---similar to standard assumptions
in random matrix theory.

In this paper, we consider the case where the data is of the type
``information${}+{}$noise.'' In other words, each observation is the sum of
two independent elements: one sampled from a ``low-dimensional''
structure, the signal part of the data, the other being
high-dimensional noise, normalized to not overwhelm but still affect
the signal. We consider two types of noise, spherical and elliptical.

In the spherical setting, we show that the spectral properties of
kernel random matrices can be understood from a new kernel matrix,
computed only from the signal part of the data, but using (in general)
a slightly different kernel. The Gaussian kernel has some special
properties in this setting.

The elliptical setting, which is important from a robustness
standpoint, is less prone to easy interpretation.
\end{abstract}

%

\end{frontmatter}

\section{Introduction}
Kernel techniques are now a standard tool of statistical practice and
kernel versions of many methods of classical multivariate statistics
have now been created. A few important examples can be found in
\citet{schoelkopf02learning} (see the description of kernel PCA, pages
41--45) and \citet{BachJordanKernelICA} (for kernel ICA), for instance.
There are several ways to describe kernel methods, but one of them is
to think of them as classical multivariate techniques using generalized
notions of inner-product. A basic input in these techniques is a kernel
matrix, that is, an inner-product (or Gram) matrix, for generalized
inner-products. If our vectors of observations are $X_1,\ldots,X_n$,
the kernel matrices studied in this paper have $(i,j)$ entry
$f(\Vert X_i-X_j \Vert_2^2)$ or $f(X_i^{\prime}X_j)$, for a certain $f$.
Popular examples include the Gaussian kernel [entries
$\exp(-\Vert X_i-X_j \Vert_2^2/2s^2)$], the Sigmoid kernel [entries
$\tanh(\kappa X_i^{\prime}X_j+\theta)$] and polynomial kernels [entries
$(X_i^{\prime}X_j)^d$]. We refer to \citet
{RasmussenWilliamsGPMLBook06} for
more examples. As explained in, for instance,
\citet{schoelkopf02learning}, kernel techniques allow
practitioners to
essentially do multivariate analysis in infinite-dimensional spaces, by
embedding the data in a infinite-dimensional space through the use of
the kernel. A nice numerical feature is that the embedding need not be
specified, and all computations can be made using the
finite-dimensional kernel matrix. Kernel techniques also allow users to
do certain forms of nonlinear data analysis and dimensionality
reduction, which is naturally very desirable.
\citet{ZwaldBousquetBlanchard04} and
\citet{vonLuxburgConsistencySpecClustering} are two interesting
relatively recent papers concerned broadly speaking with the same types
of inferential questions we have in mind and investigate in this paper,
though the settings of these papers is quite different from the one we
will work under.

Kernel matrices and the closely related Laplacian matrices also play a
central role in manifold learning [see, e.g., \citet
{BelkinNiyogiLaplacianEigenmaps03} and \citet
{IzenmanModernMultivariate08} for an overview of various techniques].
In ``classical'' statistics, they have been a mainstay of spatial
statistics and geostatistics in particular [see \citet{CressieBook93}].

In geostatistical applications, it is clear that the dimension of the
data is at most~3. Also, in applications of kernel techniques and
manifold learning, it is often assumed that the data live on a
low-dimensional manifold or structure, the kernel approach allowing us
to somehow recover (at least partially) this information. Consequently,
most theoretical analyses of kernel matrices and kernel or manifold
learning techniques have focused on situations where the data is
assumed to live on such a low-dimensional structure. In particular, it
is often the case that asymptotics are studied under the assumption
that the data is i.i.d. from a fixed distribution---independent of the
number of points. Some remarkable results have been obtained in this
setting [see \citet{KoltchinGine00} and also \citet
{BelkinNiyogiCvLapEigenmaps08}].

Let us give a brief overview of such results. In \citet
{KoltchinGine00}, the authors prove that if $X_i$ are i.i.d. with
distribution $P$, under regularity conditions on the kernel $k(x,y)$, the
$k$th largest eigenvalue of the kernel matrix $M$, with entries
\[
M(i,j)=\frac{1}{n} k(X_i,X_j) ,
\]
converges to the $k$th largest eigenvalue of the operator $K$
defined as
\[
Kf(x)=\int k(x,y)f(y)\,dP(y) .
\]
In this important paper, the authors were also able to obtain
fluctuation behavior for these eigenvalues, under certain technical
conditions [see Theorem 5.1 in \citet{KoltchinGine00}]. Similar
first-order convergence results were obtained, at a heuristic level but
through interesting arguments, in \citet{williamsseeger00}.

These results gave theoretical confirmation to practitioners' intuition
and heuristics that the kernel matrix could be used as a good proxy for
the operator $K$ on $L^2(dP)$, and hence kernel techniques could be
explained and justified through the spectral properties of this operator.

To statisticians well versed in the theory of random matrices, this set
of results appears to be similar to results for low-dimensional
covariance matrices stating that when the dimension of the data is
fixed and the number of observations goes to infinity, the sample
covariance matrix is a spectrally consistent estimator of the
population covariance matrix [see, e.g., \citet{anderson03}].
However, it is well known [see, e.g., \citet{mp67},
\citet{bai99}, \citet{JohnstoneReview07}] that this is not
the case when the dimension of the data, $p$, changes with $n$, the
number of observations, and in particular when asymptotics are studied
under the assumption that $p/n$ has a finite limit. We refer to the
asymptotic setting where $p$ and $n$ both tend to infinity as the
``high-dimensional'' setting. We note that given that more and more
datasets have observations that are high dimensional, and kernel
techniques are used on some of them [see
\citet{williamsseeger00}], it is natural to study kernel random
matrices in the high-dimensional setting.

Another important reason to study this type of asymptotics is that by
keeping track of the effect of the dimension of the data, $p$, and of
other parameters of the problem on the results, they might help us give
more accurate prediction about the finite-dimensional behavior of
certain statistics than the classical ``small $p$, large $n$''
asymptotics. An example of this phenomenon can be found in the paper
\citet{imj} where it turned out in simulation that some of the
doubly asymptotic results concerning fluctuation behavior of the
largest eigenvalue of a Wishart matrix with identity covariance are
quite accurate for $p$ and $n$ as small as 5 or 10, at least in the
right tail of the distribution. [We refer the interested reader to
\citet{imj} for more details on the specific example we just
described.] Hence, it is also potentially practically important to
carry out these theoretical studies for they can be informative even
for finite-dimensional considerations.

The properties of kernel random matrices under classical random matrix
assumptions have been studied by the author in the recent \citet
{nekKernels}. It was shown there that when the data is
high dimensional, for instance $X_i\sim\mathcal{N}(0,\Sigma_p)$, and the
operator norm of $\Sigma_p$ is, for example, bounded, kernel random
matrices essentially act like standard Gram/``covariance matrices,'' up
to recentering and rescaling, which depend only on $f$. Naturally, a
certain scaling is needed to make the problem nondegenerate, and the
results we just stated hold, for instance, when $M(i,j)=f(\Vert X_i-X_j
\Vert_2^2/p)$, for otherwise the kernel matrix is in general
degenerate. We refer to \citet{nekKernels} for more details and
discussions of the relevance of these results in practice. In limited
simulations, we found that the theory agreed with the numerics even
when $p$ was of the order of several 10's and $p/n$ was not ``too
small''
(e.g., $p/n\simeq0.2$). These results came as somewhat of a surprise
and seemed to contradict the intuition and numerous positive practical
results that have been obtained, since they suggested that the kernel
matrices we considered were just a (centered and scaled) version of the
matrix $XX^{\prime}$. However, it should be noted that the assumptions
implied that the data was truly high dimensional.

So an interesting middle ground, from modeling, theoretical and
practical points of view is the following: what happens if the data
does not live exactly on a fixed-dimensional manifold, but lives
``nearby?'' In other words, the data is now sampled from a ``noisy''
version of the manifold. This is the question we study in this paper.
We assume now that the data points $X_i \in\mathbb{R}^p$ we observe
are of the form
\[
X_i=Y_i+Z_i ,
\]
where $Y_i$ is the ``signal'' part of the observations (and live, for
instance, on a low-dimensional manifold, e.g., a three-dimensional sphere)
and $Z_i$ is the noise part of the observations (and is, e.g.,
multivariate Gaussian in dimension $p$, where $p$ might be 100).

We think this is interesting from a practical standpoint because the
assumption that the data is exactly on a manifold is perhaps a bit
optimistic and the ``noisy manifold'' version is perhaps more in line
with what statisticians expect to encounter in practice (there is a
clear analogy with linear regression here). From a theoretical
standpoint, such a model allows us to bridge the two extremes between
truly low-dimensional data and fully high-dimensional data. From a
modeling standpoint, we propose to scale the noise so that its norm
stays bounded (or does not grow too fast) in the asymptotics. That way,
the ``signal'' part of the data is likely to be affected but not totally
drowned by the noise. It is important to note, however, that the noise is
not ``small'' in any sense of the word---it is of a size comparable with
that of the signal.

In the case of spherical noise (see below for details but note that the
Gaussian distribution falls into this category) our results say that,
to first-order, the kernel matrix computed from information${}+{}$noise data
behaves like a kernel matrix computed from the ``signal'' part of the
data, but, we might have to use a different kernel than the one we
started with. This other kernel is quite explicit. In the case of
dot-product kernel matrices [i.e., $M(i,j)=f(X_i^{\prime}X_j)/n$], the
original kernel can be used (under certain assumptions)---so, to
first-order, the noise part has no effect on the spectral properties of
the kernel matrix. The results are different when looking at Euclidean
distance kernels [i.e., $M(i,j)=f(\Vert X_i-X_j \Vert_2^2)/n$] where the
effect of the noise is basically to change the kernel that is used.
This is in any case a quite positive result in that it says that the
whole body of work concerning the behavior of kernel random matrices
with low-dimensional input data can be used to also study the
``information${}+{}$noise'' case---the only change being a change of kernels.

The case of elliptical noise is more complicated. The dot-product
kernels results still have the same interpretation. But the Euclidean
distance kernels results are not as easy to interpret.


\section{Results}\label{sec2}

Before we start, we set some notation. We use $\Vert M \Vert_F$ to denote
the Frobenius norm of the matrix $M$ [so $\Vert M \Vert_F^2=\sum
_{i,j}M^2(i,j)$] and $\interleave M\interleave_2$ to denote its
operator norm, that is,
its largest singular value. We also use $\Vert\mathbf{v} \Vert_2$ to denote
the Euclidean norm of the vector $\mathbf{v}$. $a\vee b$ is shorthand
for $\max(a,b)$. Unless otherwise noted, functions that are said to be
Lipschitz are Lipschitz with respect to Euclidean norm.

We split our results into two parts, according to distributional
assumptions on the noise. One deals with the Gaussian-like case, which
allows us to give a simple proof of the results. The second part is
about the case where the noise has a distribution that satisfies
certain concentration and ellipticity properties. This is more general
and brings the geometry of the problem forward. It also allows us to
study the robustness (and lack thereof) of the results to the
sphericity of the noise, an assumption that is implicit in the
high-dimensional Gaussian (and Gaussian-like) case.

We draw some practical conclusions from our results for the case of
spherical noise in Section \ref{subsec:practicalConsequences}.

\subsection{The case of Gaussian-like noise}\label{subsec:GaussianCase}
We first study a setting where the noise is drawn according to a
distribution that is similar to a Gaussian, but slightly more general.
\begin{theorem}\label{thm:InfoPlusNoiseGaussianCase}
Suppose we observe data $X_1,\ldots,X_n$ in $\mathbb{R}^p$, with
\[
X_i=Y_i+\frac{Z_i}{\sqrt{p}} ,
\]
where $Z_i=\Sigma_p^{1/2} U_i$ where the $p$-dimensional vector $U_i$
has i.i.d. entries with mean~0, variance 1, and fourth moment
$\mu_4$, and $\{Y_i\}_{i=1}^n\sim P_n$. We assume that there
exists a deterministic vector $a$ and a real $C_1>0$, possibly
dependent on $n$, such that $\forall i, \mathbf{E}(\Vert Y_i-a \Vert
_2^2)<C_1$.
Also, $\mu_4$ might change with $n$ but is assumed to remain bounded.

$\{Z_i\}_{i=1}^n$ are i.i.d., and we also assume that $\{Y_i\}_{i=1}^n
$ and $\{Z_i\}_{i=1}^n$ are independent.

We consider the random matrices $M_f$ with $(i,j)$ entry
\[
M_f(i,j)=\frac{1}{n}f (\Vert X_i-X_j \Vert_2^2 ) \qquad \mbox
{for functions } f \in\mathcal{F}_{C_0(n)} ,
\]
where
\[
\mathcal{F}_{C_0(n)}=\Bigl\{f \mbox{ such that } {\sup_{x,y}}
|f(x)-f(y)|\leq
C_0(n)|x-y|\Bigr\}.
\]
%
Let us call $\nu=\frac{\operatorname{trace}(\Sigma_p)}{p}$. Let
$\widetilde{M}_f$
be the matrix with $(i,j)$th entry
\[
\widetilde{M}_f(i,j)=
\cases{
\dfrac{1}{n}f (\Vert Y_i-Y_j \Vert_2^2+2\nu), &\quad if
$i\neq j$,\vspace*{1pt}\cr
\dfrac{1}{n} f(0), &\quad if $i=j$.}
\]
Assuming only that $\mu_4$ is bounded uniformly in $n$, we
have, for a constant $C$ independent of $n$, $p$ and $\Sigma_p$,
%
%
\begin{equation}\label{eq:ExplicitControlErrorGaussianLikeCase}
\mathbf{E}^*\Bigl({\sup_{f \in\mathcal{F}_{C_0(n)}}}\Vert M_f-\widetilde
{M} _f \Vert_F^2\Bigr)\leq
C C_0^2(n) \biggl[ \frac{\operatorname{trace}(\Sigma
_p^2)}{p^2}+\frac{\interleave\Sigma _p\interleave_2}{p} C_1 \biggr] .
\end{equation}

We place ourselves in the high-dimensional setting where $n$ and $p$
tend to infinity.
We assume that $\operatorname{trace}(\Sigma_p^2)/p^2\rightarrow0$,
as $p$ tends to
infinity.

Under these assumptions, for any fixed $C_0>0$ and $C_1>0$,
\[
{\lim_{n,p\rightarrow\infty}} \sup_{f\in\mathcal{F}_{C_0}}\Vert
M_f-\widetilde{M} _f \Vert_F^2=0 \qquad\mbox{in probability}.
\]
\end{theorem}

If we further assume that $\nu$ remains, for instance, bounded,
the same result holds if we replace the
diagonal of $\widetilde{M}$ by $f(2\nu)/n$, because
$|f(2\nu)-f(0)|\leq2\nu C_0$ and therefore ${\sup
_{f\in\mathcal{F}_{C_0}}} |f(2\nu)-f(0)|\leq2\nu C_0$.
The approximating matrix we then get is the matrix with $(i,j)$th entry
$f_{\nu}(\Vert Y_i-Y_j \Vert_2^2)$, where $f_{\nu
}(x)=f(x+2\nu)$, that is, a ``pure signal'' matrix
involving a different kernel from the one with which we started.

We note that there is a potential measurability issue that we address
in the proof. Our theorem really means that we can find a random
variable that dominates the ``random element'' ${\sup_{f\in\mathcal
{F}_{C_0(n)}}}\Vert M_f-\widetilde{M}_f \Vert_F^2$ and goes to 0 in
probability.
(This measurability issue could also be addressed through separability
arguments but outer-probability statements suffice for our purposes in
this paper.)

A subcase of our result is the case of Gaussian noise: then $U_i$ is
$\mathcal{N}(0,\operatorname{Id}_p)$ and our result naturally applies.

We also note that $P_n$ can change with $n$. The class of functions we
consider is fixed in the last statement of
the theorem but if we were to look at a sequence of kernels we could
pick a different function in the class $\mathcal{F}_{C_0}$ for each $n$
[the proof also applies to matrices with entries $M(i,j)=f_{(i,j)}
(\Vert X_i-X_j \Vert_2^2 )/n$, where the functions considered also
depend on $(i,j)$, but we present the results with a function $f$
common to all entries]. It should also be noted that the proof
technique allows us to deal with classes of functions that vary with
$n$: we could have a varying $C_0(n)$. As (\ref
{eq:ExplicitControlErrorGaussianLikeCase}) makes clear, the
approximation result will hold
as soon as the right-hand side of (\ref
{eq:ExplicitControlErrorGaussianLikeCase}) goes to 0 asymptotically,
that is,
$C_0^2(n)\max(\operatorname{trace}(\Sigma_p^2)/p^2,\interleave
\Sigma _p\interleave_2/p)\rightarrow0$.
Finally, we work here with uniformly Lipschitz functions. The proof
technique carries over to other classes, such as certain classes of
H\"{o}lder functions, but the bounds would be different.
\begin{pf*}{Proof of Theorem \ref{thm:InfoPlusNoiseGaussianCase}}
The strategy is to use the same entry-wise expansion approach that\vspace*{1pt} was
used in \citet{nekKernels}. To do so, we remark that $\Vert
Z_i-Z_j \Vert_2^2/p$ remains essentially constant [across $(i,j)$] in the
setting we are considering---this is a consequence of the ``spherical''
nature of high-dimensional Gaussian distributions. We can therefore try
to approximate $M(i,j)$ by $f(\Vert Y_i-Y_j \Vert_2^2+2\nu)/n$ and
all we need to do is to show that the remainder is small.

We also note that if, as we assume, $\operatorname{trace}(\Sigma
_p^2)/p^2\rightarrow0$,
then $\interleave\Sigma_p\interleave_2=o(p)$, since $\interleave
\Sigma_p\interleave_2^2\leq
\operatorname{trace}(\Sigma_p^2)$.

- \textit{Work conditional on} $\mathcal{Y}_n=\{
Y_i\}
_{i=1}^n$, \textit{for} $i\neq j$.

We clearly have
\[
\Vert X_i-X_j \Vert_2^2=\Vert Y_i-Y_j \Vert_2^2+2\frac
{(Z_i-Z_j)^{\prime}}{\sqrt{p}}
(Y_i-Y_j)+\frac{\Vert Z_i-Z_j \Vert_2^2}{p} .
\]
Let us study the various parts of this expansion. Conditional on
$\mathcal{Y}_n$, if we call $y_{i,j}=Y_i-Y_j$, we see easily that
\[
Z_i-Z_j=\Sigma_p^{1/2}(U_i-U_j)
\]
and
\[
(Z_i-Z_j)^{\prime}(Y_i-Y_j)=(U_i-U_j)^{\prime}\Sigma
_p^{1/2}y_{i,j}.
\]
Note that $U_i-U_j$, which we denote $\Gamma_{i,j}$, has i.i.d.
entries, with mean 0, variance~$2$ and fourth moment $2\mu_4+6$.
We call
\[
\alpha_{i,j}=(Z_i-Z_j)^{\prime}(Y_i-Y_j)/\sqrt{p}
\]
and
\[
\beta
_{i,j}=\frac{\Vert Z_i-Z_j \Vert_2^2}{p}-2\frac{\operatorname
{trace}(\Sigma_p)}{p}.
\]

With this notation, we have
\[
\Vert X_i-X_j \Vert_2^2-(\Vert Y_i-Y_j \Vert_2^2+2\nu)=2\alpha
_{i,j}+\beta_{i,j} .
\]

Therefore, for any function $f$ in $\mathcal{F}_{C_0(n)}$,
%
\[
\bigl|f(\Vert X_i-X_j \Vert_2^2)-f(\Vert Y_i-Y_j \Vert_2^2+2\nu)\bigr|\leq
C_0(n) |\beta_{i,j}+2\alpha_{i,j} | ,
\]
and hence,
\[
[f(\Vert X_i-X_j \Vert_2^2)-f(\Vert Y_i-Y_j \Vert_2^2+2\nu
)
]^2 \leq2 C_0(n)^2 [\beta_{i,j}^2+4\alpha_{i,j}^2 ] .
\]
We naturally also have
\[
\sup_{f\in\mathcal{F}_{C_0(n)}} [f(\Vert X_i-X_j \Vert
_2^2)-f(\Vert Y_i-Y_j \Vert_2^2+2\nu) ]^2 \leq2 C_0(n)^2
[\beta
_{i,j}^2+4\alpha_{i,j}^2 ] .
\]
So we have found a random variable $\tau_n=2 C_0^2(n) [\beta
_{i,j}^2+4\alpha_{i,j}^2 ]$ that dominates the random element
$\zeta_n=\sup_{f\in\mathcal{F}_{C_0(n)}} [f(\Vert X_i-X_j
\Vert_2^2)-f(\Vert Y_i-Y_j \Vert_2^2+2\nu) ]^2$. One
might be concerned about the measurability of $\zeta_n$---but by using
outer expectations [see \citet{vandervaart}, page 258], we can
completely bypass this potential problem. In what follows, we denote by
$\mathbf{E}^*(\cdot)$ an outer expectation. (Though this technical point
does not shed further light on the problem, it naturally needs to be addressed.)

Hence,
\begin{eqnarray*}
&&\mathbf{E}^*\Bigl(\sup_{f\in\mathcal{F}_{C_0(n)}} \bigl(f(\Vert X_i-X_j
\Vert_2^2)-f(\Vert Y_i-Y_j \Vert_2^2+2\nu) \bigr)^2|\mathcal
{Y}_n\Bigr) \\
&&\qquad\leq2 C_0(n)^2 \bigl(\mathbf{E}(\beta_{i,j}^2)+\mathbf{E}(4\alpha
_{i,j}^2|\mathcal{Y}_n) \bigr) .
\end{eqnarray*}

Let us focus on $\mathbf{E}(\beta_{i,j}^2)$ for a moment. Let us call
$\Gamma
_{i,j}=U_i-U_j$. We first note that $\Vert Z_i-Z_j \Vert_2^2=\Gamma
_{i,j}^{\prime}\Sigma_p \Gamma_{i,j}=\operatorname{trace}(\Sigma_p
\Gamma_{i,j}\Gamma_{i,j}^{\prime})$. In particular,
\[
\mathbf{E}(\Vert Z_i-Z_j \Vert_2^2)=2\operatorname{trace}(\Sigma
_p) ,
\]
so $\mathbf{E}(\beta_{i,j})=0$. Therefore, $\mathbf{E}(\beta
_{i,j}^2)=\operatorname{var}(\Vert Z_i-Z_j \Vert_2^2)/p^2$. Now
recall the results found, for
instance, in Lemma
A-1 in \citet{nekKernels}: if the vector $\gamma$ has i.i.d.
entries with mean 0, variance $\sigma^2$ and fourth moment $\kappa_4$,
and if $M$ is a symmetric matrix,
\[
\mathbf{E}((\gamma^{\prime}M \gamma)^2)=\sigma^4\bigl(2\operatorname
{trace}(M^2)+\operatorname{trace}(M)^2\bigr)+(\kappa_4-3\sigma
^4)\operatorname{trace}(M\circ M) ,
\]
where $M\circ M$ is the Hadamard product of $M$ with itself, that is,
the entrywise product of two matrices.

Applying this result in our setting [i.e., using the moments (given
above) of $\Gamma_{i,j}$, which has i.i.d. entries, in the previous
formula] gives
\[
\operatorname{var}(\Vert Z_i-Z_j \Vert_2^2)=\operatorname
{var}(\Gamma_{i,j}^{\prime}\Sigma_p \Gamma_{i,j})=8 \operatorname
{trace}(\Sigma_p^2)+2(\mu_4-3)\operatorname{trace}(\Sigma_p\circ
\Sigma_p) .
\]
It is easy to see that $\operatorname{trace}(\Sigma_p\circ\Sigma
_p)\leq\operatorname{trace}(\Sigma_p^2)$, since $\operatorname
{trace}(\Sigma_p^2)=\sum_{i,j}\sigma_p^2(i,j)$ and
$\operatorname{trace}(\Sigma_p\circ\Sigma_p)=\sum_{i}\sigma
_p^2(i,i)$. Therefore,
\[
\mathbf{E}(\beta_{i,j}^2)=\frac{\operatorname{var}(\Vert Z_i-Z_j
\Vert_2^2)}{p^2}\leq
\frac{8+2(\mu_4-3)}{p^2} \operatorname{trace}(\Sigma_p^2)=O
\biggl(\frac
{\operatorname{trace}(\Sigma_p^2)}{p^2} \biggr) .
\]
We note that under our assumptions on $\operatorname{trace}(\Sigma
_p^2)/p^2$ and the
fact that $\mu_4$ remains bounded in $n$ (and therefore $p$),
this term will go to 0 as $p\rightarrow\infty$.

On the other hand, because $\alpha_{i,j}|\mathcal{Y}_n=\Gamma
_{i,j}^{\prime}
\Sigma_p^{1/2} y_{i,j}/\sqrt{p}$, and because $\mathbf{E}(\Gamma_{i,j})=0$
and $\operatorname{cov}(\Gamma_{i,j})=2\operatorname{Id}_p$, we have
\[
\mathbf{E}(\alpha_{i,j}^2|\mathcal{Y}_n)=2 \frac{y_{i,j}^{\prime
}\Sigma_p
y_{i,j}}{p}\leq2\interleave\Sigma_p\interleave_2 \frac{\Vert
y_{i,j} \Vert
_2^2}{p}\leq4
\interleave\Sigma_p\interleave_2\frac{\Vert Y_i-a \Vert
_2^2+\Vert Y_j-a \Vert
_2^2}{p} .
\]

Hence, we have for $C$ a constant independent of $\Sigma_p$, $p$ and $n$,
\begin{eqnarray*}
&&\mathbf{E}^*\Bigl(\sup_{f\in\mathcal{F}_{C_0(n)}} \bigl(f(\Vert X_i-X_j
\Vert_2^2)-f(\Vert Y_i-Y_j \Vert_2^2+2\nu) \bigr)^2|\mathcal
{Y}_n\Bigr) \\
&&\qquad\leq C C_0^2(n) \biggl[ \frac{\operatorname{trace}(\Sigma
_p^2)}{p^2}+\frac{\interleave\Sigma_p\interleave_2}{p} (\Vert
Y_i-a \Vert_2^2+\Vert Y_j-a \Vert
_2^2
) \biggr] .
\end{eqnarray*}
This inequality allows us to conclude that, for another constant $C$,
\[
\mathbf{E}^*\Bigl({\sup_{f\in\mathcal{F}_{C_0(n)}}} \Vert M_f-\widetilde
{M} _f \Vert_F^2|\mathcal{Y}_n\Bigr) \leq C C_0^2(n) \Biggl[ \frac
{\operatorname{trace}(\Sigma_p^2)}{p^2}+\frac{\interleave\Sigma
_p\interleave_2}{p} \frac{1}{n}\sum_{i=1}^n \Vert Y_i-a \Vert_2^2
\Biggr] ,
\]
since clearly,
\[
{\sup_{f\in\mathcal{F}_{C_0(n)}}} \Vert M_f-\widetilde{M}_f \Vert
_F^2 \leq
\frac
{1}{n^2}\sum_{i,j} \sup_{f\in\mathcal{F}_{C_0(n)}} \bigl(f(\Vert
X_i-X_j \Vert_2^2)-f(\Vert Y_i-Y_j \Vert_2^2+2\nu) \bigr)^2 .
\]

Under the assumption that $\mathbf{E}(\Vert Y_i-a \Vert_2^2)$ exists
and is less
than $C_1$, we finally conclude that
\[
\mathbf{E}^*\Bigl({\sup_{f \in\mathcal{F}_{C_0(n)}}}\Vert M_f-\widetilde
{M} _f \Vert_F^2\Bigr)\leq C
C_0^2(n) \biggl[ \frac{\operatorname{trace}(\Sigma_p^2)}{p^2}+\frac
{\interleave\Sigma _p\interleave_2}{p} C_1 \biggr],
\]
and (\ref{eq:ExplicitControlErrorGaussianLikeCase}) is shown.

Therefore, under our assumptions,
\[
\mathbf{E}^*\Bigl({\sup_{f \in\mathcal{F}_{C_0}}}\Vert M_f-\widetilde{M}
_f \Vert_F^2\Bigr)=o(1) .
\]
Hence, when $n$ and $p$ tend to $\infty$,
\[
\sup_{f \in\mathcal{F}_{C_0}} \Vert M-\widetilde{M} \Vert_F^2
\rightarrow0\qquad
\mbox{in
probability,}
\]
as announced in the theorem.
\end{pf*}

\subsection{Case of noise drawn from a distribution satisfying
concentration inequalities}\label{subsec:ConcCase}

The proof of Theorem \ref{thm:InfoPlusNoiseGaussianCase} makes clear
that the heart of our argument is geometric: we exploit the fact that
$\Vert Z_i-Z_j \Vert_2^2/p$ is essentially constant across pairs
$(i,j)$. It
is therefore natural to try to extend the theorem to more general
assumptions about the noise distribution than the Gaussian-like one we
worked under previously. It is also important to understand the impact
of the implicit geometric assumptions (i.e., sphericity of the noise)
that are made and in particular the robustness of our results against
these geometric assumptions.

We extend the results in two directions. First, we investigate the
generalization of our Gaussian-like results to the setting of
Euclidean-distance kernel random matrices, when the noise is
distributed according to a distribution satisfying a concentration
inequality multiplied by a random variable, that is, a generalization
of elliptical distributions. This allows us to show that the
Gaussian-like results of Theorem \ref{thm:InfoPlusNoiseGaussianCase}
essentially hold under much weaker assumptions on the noise
distribution, as long as the Gaussian geometry (i.e., a spherical
geometry) is preserved (see Corollary \ref
{coro:SphericalNoiseKernels}). The results of Theorem \ref
{thm:InfoPlusNoiseConcCase} show that breaking the Gaussian geometry
results in quite different approximation results.

We also discuss in Theorem
\ref{thm:InfoPlusNoiseConcCaseDotProductsKernels} the situation of
inner-product kernel random matrices under the same ``generalized
elliptical'' assumptions on the noise.

\subsubsection{The case of Euclidean distance kernel random matrices}
We have the following theorem.
\begin{theorem}[(Euclidean distance kernels)]\label{thm:InfoPlusNoiseConcCase}
Suppose we observe data $X_1,\ldots, X_n$ in $\mathbb{R}^p$, with
\[
X_i=Y_i+R_i \frac{Z_i}{\sqrt{p}} .
\]
We place ourselves in the high-dimensional setting where $n$ and $p$
tend to infinity.
We assume that $\{Y_i\}_{i=1}^n\sim P_n$.

$\{Z_i\}_{i=1}^n$ are i.i.d. with $\mathbf{E}(Z_i)=0$, and we also assume
that $\mathcal{Y}_n=\{Y_i\}_{i=1}^n $ and $\{Z_i\}_{i=1}^n$ are
independent. $R_i$ are random variables independent of $\{Z_i\}_{i=1}^n$.

We now assume that the distribution of $Z_i$ is such that, for any
1-Lipschitz function $F$, if $\mu_F=\mathbf{E}(F(Z_i))$,
\[
P\bigl(|F(Z_i)-\mu_F|>r\bigr)\leq C \exp(-c_0 r^b)\triangleq h(r) ,
\]
where for simplicity we assume that $c_0$, $C$ and $b$ are independent
of $p$. We call $\nu=\mathbf{E}(\Vert Z_i \Vert_2^2)/p$ and assume that
$\nu$ stays bounded as $p\rightarrow\infty$.

We assume that $\forall i$, $|R_i|\in[r_{\infty}(p),R_{\infty}(p)]$,
where $r_{\infty}(p)$ and $R_{\infty}(p)$ are deterministic sequences
depending on $p$. We assume without loss of generality that $R_{\infty
}(p)\geq1$.

Calling $\mathcal{M}(\mathcal{Y}_n)=\max_{i\neq j} \Vert Y_i-Y_j
\Vert_2^2$, we
assume that there exists $\mathcal{M}_p$ such that $P(\mathcal
{M}(\mathcal{Y}_n)\leq\mathcal{M}_p)\rightarrow1$ and $\varepsilon>0$
such that
\[
\max(\mathcal{M}_p^{1/2},R_{\infty}(p)) \frac{R_{\infty}(p) (\log
n+(\log
n)^{\varepsilon})^{1/b}}{\sqrt{p}}\rightarrow0.
\]

Then we have
%
%
\begin{equation}\label{eq:InterpointDistanceEllipCase}\quad
\max_{i\neq j} \bigl|\Vert X_i-X_j \Vert_2^2- [\Vert Y_i-Y_j
\Vert_2^2+\nu(R_i^2+R_j^2) ] \bigr|\rightarrow0\qquad
\mbox{in probability}.
\end{equation}

We call $\mathcal{W}(\mathcal{Y}_n)={\min_{i\neq j}}\Vert Y_i-Y_j
\Vert_2^2$, and
suppose we pick $\mathcal{W}_p$ such that $P(\mathcal{W}(\mathcal
{Y}_n)\geq
\mathcal{W}_p)\rightarrow1$. (Note that $\mathcal{W}_p=0$ is always a
possibility.)

We call, for $\eta>0$ given, $I_p(\eta)=[\mathcal
{W}_p+2\nu
r_{\infty}^2(p)-\eta,\mathcal{M}_p+2\nu R_{\infty
}^2(p)+\eta]$, and
\[
\mathcal{F}_{C_1,I_p(\eta)}= \Bigl\{f \mbox{ such that } {\sup_{x,y
\in
I_p(\eta)}} |f(x)-f(y)| \leq C_1|x-y| \Bigr\}.
\]
We consider the random matrices $M_f$ with $(i,j)$ entry
\[
M_f(i,j)=\frac{1}{n}f (\Vert X_i-X_j \Vert_2^2 ) \qquad \mbox{for
} f
\in\mathcal{F}_{C_1,I_p(\eta)} .
\]

Let us call $\widetilde{M}_f$ the matrix with $(i,j)$th entry
\[
\widetilde{M}_f(i,j)=
\cases{
\dfrac{1}{n}f \bigl(\Vert Y_i-Y_j \Vert_2^2+\nu
(R_i^2+R_j^2) \bigr),
&\quad if $i\neq j$,\vspace*{1pt}\cr
\dfrac{1}{n} f(0), &\quad if $i=j$.}
\]
We have, for any given $C_1>0$ and $\eta>0$,
%
%
\begin{equation}\label{eq:MainResConcEllipCase}
{\lim_{n,p\rightarrow\infty}} \sup_{f\in\mathcal{F}_{C_1,I_p(\eta)}}
\Vert M_f-\widetilde{M}_f \Vert_F=0 \qquad\mbox{in probability}.
\end{equation}
\end{theorem}

We have the following corollary in the case of ``spherical'' noise,
which is a generalization of the Gaussian-like case considered in
Theorem \ref{thm:InfoPlusNoiseGaussianCase}.
\begin{corollaryNumThm}[(Euclidean distance kernels with spherical
noise)]\label{coro:SphericalNoiseKernels}
Suppose we observe data $X_1,\ldots,X_n$ in $\mathbb{R}^p$, with
\[
X_i=Y_i+\frac{Z_i}{\sqrt{p}} ,
\]
where $Y_i$ and $Z_i$ satisfy the same assumptions as in
Theorem \ref{thm:InfoPlusNoiseConcCase} [with $r_{\infty}(p)=R_{\infty}(p)=1$].
Then the results of Theorem \ref{thm:InfoPlusNoiseConcCase} apply with
\[
I_p(\eta)=[\mathcal{W}_p+2\nu-\eta,\mathcal
{M}_p+2\nu
+\eta]
\]
and
\[
\widetilde{M}_f(i,j)=
\cases{
\dfrac{1}{n}f (\Vert Y_i-Y_j \Vert_2^2+2\nu), &\quad if
$i\neq j$,\vspace*{2pt}\cr
\dfrac{1}{n} f(0), &\quad if $i=j$.}
\]
\end{corollaryNumThm}

As in Theorem \ref{thm:InfoPlusNoiseGaussianCase}, we deal with
potential measurability issues concerning the $\sup$ in the proof. Our
theorem is really that we can find a random variable that goes to 0
with probability 1 and dominates\vspace*{-2pt} the random element ${\sup_{f\in
\mathcal{F}_{C_1,I_p(\eta)}} }\Vert M_f-\widetilde{M}_f \Vert_F$---an
outer-probability statement.

This theorem generalizes Theorem \ref{thm:InfoPlusNoiseGaussianCase} in
two ways. The ``spherical'' case, detailed in Corollary \ref
{coro:SphericalNoiseKernels}, is a more general version of Theorem
\ref{thm:InfoPlusNoiseGaussianCase} limited to Gaussian noise. This is
because the Gaussian setting corresponds to $b=2$ and
$c_0=1/(2\interleave\Sigma_p\interleave_2)$. However, assuming
``only'' concentration inequalities
allows us to handle much more complicated structures for the noise
distribution. Some examples are given below. We also note that if the
$Y_i$'s (i.e., the signal part of the $X_i$'s) are sampled, for
instance, from a fixed manifold of finite Euclidean diameter, the
conditions on $\mathcal{M}$ are automatically satisfied, with
$\mathcal{M}_p$
being the Euclidean diameter of the corresponding manifold.

Another generalization is ``geometric'': by allowing $R_i$ to vary with
$i$, we move away from the spherical geometry of high-dimensional
Gaussian vectors (and generalizations), to a more ``elliptical''
setting. Hence, our results show clearly the potential limitations and
the structural assumptions that are made when one assumes Gaussianity
of the noise. Theorem \ref{thm:InfoPlusNoiseConcCase} and Corollary
\ref
{coro:SphericalNoiseKernels} show that the Gaussian-like results of
Theorem \ref{thm:InfoPlusNoiseGaussianCase} are not robust against a
change in the geometry of the noise. We note however that if $R_i$ is
independent of $Z_i$ and $\mathbf{E}(R_i^2)=1$, $\operatorname
{cov}(R_i Z_i)=\operatorname{cov}(Z_i)$,
so all the noise models have the same covariance but they may yield
different approximating matrices and hence different spectral behavior
for our information${}+{}$noise models.

However, the spherical results have the advantage of having simple
interpretations. In the setting of Corollary \ref
{coro:SphericalNoiseKernels}, if we assume that $f(0)$ and
$f(2\nu)$ are uniformly bounded (in $n$) over the class of
functions we consider, we can replace the diagonal\vspace*{1pt} of $\widetilde{M}$ by
$f(2\nu)/n$ and have the same approximation results. Then the
``new'' $\widetilde{M}$ is a kernel matrix computed from the signal
part of
the data with the new kernel $f_{\nu}(x)=f(x+2\nu)$.

To make our result more concrete, we give a few examples of
distributions for which the concentration assumptions on $Z_i$ are satisfied:
\begin{itemize}
\item Gaussian random variables, for which we have $c_0=1/(2\interleave
\Sigma\interleave_2)$. We refer to
Ledoux [(\citeyear{ledoux2001}), Theorem 2.7] for a justification of this claim.
\item Vectors of the type $\sqrt{p} \mathbf{v}$ where $\mathbf{v}$
is uniformly
distributed on the unit ($\ell_2$-)sphere in dimension $p$. Theorem
2.3 in \citet{ledoux2001} shows that
our assumptions are satisfied, with $c(p)=(1-1/p)/2\geq c_0=1/4$, after
noticing that a 1-Lipschitz function with respect to Euclidean norm is
also 1-Lipschitz with respect
to the geodesic distance on the sphere.
\item Vectors $\Gamma\sqrt{p} \mathbf{v}$, with $\mathbf{v}$ uniformly
distributed on the unit ($\ell_2$-)sphere in $\mathbb{R}^p$
and with $\Gamma\Gamma'=\Sigma$ having bounded operator norm.
\item Vectors of the type $p^{1/b} \mathbf{v}$, $1\leq b \leq2$, where
$\mathbf{v}$ is uniformly
distributed in the unit $\ell^b$ ball or sphere in $\mathbb{R}^p$. (See
Ledoux [(\citeyear{ledoux2001}), Theorem 4.21] which refers to
\citet{SchechtmanZinn00} as the source of the theorem.) In this
case, $c_0$ depends only on $b$.
\item Vectors with log-concave density of the type $e^{-U(x)}$, with
the Hessian of $U$
satisfying, for all $x$, $\operatorname{Hess}(U)\geq2c_0\operatorname
{Id}_p$, where
$c_0>0$ is the real that appears in our assumptions. See
Ledoux [(\citeyear{ledoux2001}), Theorem 2.7] for a justification.
\item Vectors $\mathbf{v}$ distributed according to a (centered)
Gaussian copula, with corresponding correlation matrix, $\Sigma$,
having $\interleave\Sigma\interleave_2$ bounded.
We refer to \citet{nekCorrEllipD} for a justification of the fact
that our assumptions are satisfied. [If $\widetilde{\mathbf{v}}$ has a
Gaussian copula distribution,
then its $i$th entry satisfy $\widetilde{\mathbf{v}}_i=\Phi(N_i)$, where
$N$ is multivariate normal with covariance matrix $\Sigma$, $\Sigma$
being a correlation matrix,
that is, its diagonal is 1. Here $\Phi$ is the cumulative distribution
function of a standard normal distribution. Taking
$\mathbf{v}=\widetilde{\mathbf{v}}-1/2$ gives a centered Gaussian
copula.]
This last example is intended to show that the result can handle quite
complicated and nonlinear noise structure.
\end{itemize}
We note that to justify that the assumptions of the theorem are
satisfied, it is enough to be able to show concentration around the
mean or the median, as Proposition~1.8 in \citet{ledoux2001}
makes clear.

The reader might feel that the assumptions concerning the boundedness
of the $R_i$'s will be limiting in practice. We note that the same
proof essentially goes through if we just require that $|R_i|$'s belong
to the interval $[r_{\infty}(p),R_{\infty}(p)]$ with probability going
to 1, but this requires a little bit more conditioning and we leave the
details, which are not difficult, to the interested reader. So for
instance, if we had a tail condition on $|R_i|$, we could bound $\max
|R_i|$ with high probability to get a choice of $R_{\infty}(p)$. So
this boundedness condition is here just to make the exposition simpler
and is not particularly limiting in our opinion. On the other hand, we
note that our conditions allow dependence in the $R_i$'s and are
therefore rather weak requirements.

Finally, the theorem as stated is for a fixed $C_1$, though the class
of functions we are considering might vary with $n$ and $p$ through the
influence of $I_p(\eta)$. The proof makes clear that $C_1$ could also
vary with $n$ and $p$. We discuss in more details the necessary
adjustments after the proof.

\begin{pf*}{Proof of Theorem \ref{thm:InfoPlusNoiseConcCase}}
We use the notation $\mathcal{Y}_n=\{Y_i\}_{i=1}^n$ and $P_{\mathcal{Y}_n}$
to denote probability conditional on $\mathcal{Y}_n$. We call
$\mathcal{L}=\{
\mathcal{Y}_n \dvtx \mathcal{M}(\mathcal{Y}_n)\leq\mathcal{M}_p\}$.

Let us also call $\mathcal{YR}_n=\{\{Y_i\}_{i=1}^n,\{R_i\}_{i=1}^n\}$;
similarly, $P_{\mathcal{YR}_n}$ denotes probability conditional on
$\mathcal{YR}_n$. We call $\mathcal{LR}=\{\mathcal{YR}_n\dvtx \mathcal
{Y} \in\mathcal{L}\}$.
We will start by working conditionally on $\mathcal{YR}_n$ and eventually
decondition our results.

We assume from now on that the $\mathcal{YR}_n$ we work with is such that
$\mathcal{Y}_n\in\mathcal{L}$. Note that $P(\mathcal{Y}_n\in
\mathcal{L})\rightarrow
1$ by assumption and also $P(\mathcal{YR}_n \in\mathcal
{LR})\rightarrow1$.


The main idea now is that, in a strong sense,
\[
\forall i\neq j\qquad \Vert X_i-X_j \Vert_2^2\simeq\Vert Y_i-Y_j \Vert_2^2
+(R_i^2+R_j^2) \nu,
\]
where $\nu=\mathbf{E}(Z_i^2)$.
To show this formally, we write
\[
\Vert X_i-X_j \Vert_2^2- [\Vert Y_i-Y_j \Vert_2^2 +(R_i^2+R_j^2)
\nu
]=2\alpha_{i,j}+\beta_{i,j} ,
\]
where
\[
\alpha_{i,j}=\frac{(R_i Z_i-R_j Z_j)^{\prime}(Y_i-Y_j)}{\sqrt{p}}
\]
and
\[
\beta_{i,j}=\frac{\Vert R_i Z_i-R_j Z_j \Vert_2^2}{p}-\frac{(R_i^2
\mathbf{E}(\Vert Z_i \Vert_2^2)+R_j^2 \mathbf{E}(\Vert Z_j \Vert
_2^2))}{p} .
\]

Our aim is to show that, as $n$ and $p$ tend to infinity,
\[
{\max_{i\neq j}} |\alpha_{i,j}|+|\beta_{i,j}|\rightarrow0 \qquad\mbox{in
probability.}
\]

- \textit{On} ${\max_{i\neq j}} |\alpha_{i,j}|$.

Note that if $i=j$, $\alpha_{i,j}=0$. Clearly,
\begin{eqnarray*}
P_{\mathcal{YR}_n}(|\alpha_{i,j}|>2r)&\leq& P_{\mathcal{YR}_n}
\biggl(|R_i|\frac
{|Z_i^{\prime}(Y_i-Y_j)|}{\sqrt{p}}>r \biggr)\\
&&{} + P_{\mathcal
{YR}_n}
\biggl(|R_j|\frac
{|Z_j^{\prime}(Y_i-Y_j)|}{\sqrt{p}}>r \biggr).
\end{eqnarray*}
Since we assumed that $|R_i|\leq R_{\infty}(p)$, we see that the
function $F_{i,j}(Z)=R_i Z^{\prime}(Y_i-Y_j)/\sqrt{p}$ is Lipschitz (with
respect to Euclidean norm), with Lipschitz constant smaller than
$(\mathcal{M}_p)^{1/2} R_{\infty}(p)/\sqrt{p}$, when $\mathcal
{Y}_n$ is in
$\mathcal{L}$. Also, since $\mathbf{E}(Z_i)=0$, $\mathbf
{E}(F_{i,j}(Z)|\mathcal{YR}_n)=0$, where the expectation is conditional
on $\mathcal{YR}_n$.
Hence, our concentration assumptions on $Z_i$ imply that
\[
P_{\mathcal{YR}_n}\bigl(|R_i|\bigl|Z_i^{\prime}(Y_i-Y_j)/\sqrt{p}\bigr|>r\bigr)\leq
C\exp\bigl(-c_0
\bigl(p^{1/2}r/[\mathcal{M}_p^{1/2} R_{\infty}(p)]\bigr)^b\bigr) .
\]
Therefore, if we use a simple union bound, we get
\[
P_{\mathcal{YR}_n}\Bigl({\max_{i\neq j}} |\alpha_{i,j}|>2r\Bigr)\leq2C n^2 \exp\bigl(-c_0
\bigl(p^{1/2}r/[\mathcal{M}_p^{1/2} R_{\infty}(p)]\bigr)^b\bigr) .
\]
In particular, if we pick, for $\varepsilon>0$, $r_0=R_{\infty
}(p)\mathcal
{M}_p^{1/2}p^{-1/2}(\log n +(\log n)^{\varepsilon
})^{1/b}(2/\break c_0)^{1/b}$, we
see that
\begin{eqnarray*}
P_{\mathcal{YR}_n}\Bigl({\max_{i\neq j}} |\alpha_{i,j}|>2r_0\Bigr)&\leq& 2C n^2
\exp
\bigl(-c_0 \bigl(p^{1/2}r_0/[\mathcal{M}_p^{1/2}R_{\infty}(p)]\bigr)^b\bigr)\\
&=& 2C \exp
(-2(\log
n)^{\varepsilon})\rightarrow0.
\end{eqnarray*}

Since
\[
P\Bigl({\max_{i,j}} |\alpha_{i,j}| >t\Bigr)\leq P\Bigl({\max_{i,j}} |\alpha_{i,j}| >t
\mbox{ and } \mathcal{YR}_n\in\mathcal{LR}\Bigr)+P(\mathcal{YR}_n\notin\mathcal
{LR}) ,
\]
and since the latter goes to 0, we have, unconditionally,
\[
P\Bigl({\max_{i,j}} |\alpha_{i,j}| >2r_0\Bigr)\rightarrow0 .
\]

- \textit{On} ${\max_{i\neq j}} |\beta_{i,j}|$.

We see that if $A$ and $B$ are vectors in $\mathbb{R}^p$, the map
$N_{R_i,R_j}\dvtx (A,B)\rightarrow\Vert R_i A -R_j B \Vert_2$ is
$(|R_i|\vee
|R_j|)$-Lipschitz
on $\mathbb{R}^{2p}$ equipped with the norm $\Vert A \Vert_2+\Vert B
\Vert_2$, by
the triangle inequality. Therefore, using Propositions 1.11 and 1.7 in
\citet{ledoux2001} [and using the fact that $h(r)\rightarrow0$ as
$r\rightarrow\infty$ and $h$ is continuous when using the latter], we
conclude that
%
%
\begin{equation}\label{eq:ConditionalConcentration}\qquad
P_{\mathcal{YR}_n} \bigl( | \Vert R_i Z_i-R_j Z_j \Vert
_2-\mathbf{E}(\Vert R_i Z_i-R_j Z_j \Vert_2) |>r \bigr)\leq4
h\bigl(r/(2R_{\infty}(p))\bigr) .
\end{equation}
If now\vspace*{2pt} $\gamma_{i,j}=\mathbf{E}(\Vert R_i Z_i-R_j Z_j \Vert
_2|\mathcal{YR}_n)$, and
if $r_1=2R_{\infty}(p)(2/c_0)^{1/b}(\log n+(\log n)^{\varepsilon
})^{1/b}p^{-1/2}$,
\[
P_{\mathcal{YR}_n} \biggl({\max_{i\neq j}} \biggl|\frac{\Vert R_i
Z_i-R_j Z_j \Vert_2-\gamma_{i,j}}{\sqrt{p}} \biggr|>r_1 \biggr)\leq
K \exp(-(\log
n)^{\varepsilon})\rightarrow0 ,
\]
where $K$ is a constant which does not depend on $\mathcal{YR}_n$. So we
conclude that unconditionally, if
\begin{eqnarray*}
\Delta_0&=&{\max_{i\neq j} }\biggl|\frac{\Vert R_i Z_i-R_j Z_j \Vert
_2-\gamma
_{i,j}}{\sqrt{p}} \biggr| ,
\\
P (\Delta_0>r_1 )&\rightarrow&0.
\end{eqnarray*}
Note also that under our assumptions, $r_1\rightarrow0$.
Recall that we aim to show that
\[
\Delta_2={\max_{i\neq j}} \biggl|\frac{\Vert R_i Z_i-R_j Z_j \Vert
_2^2}{p}-\nu(R_i^2+R_j^2) \biggr|\rightarrow0 \qquad\mbox{in
probability.}
\]
Let us first work on
\[
\Delta_1={\max_{i\neq j}} \biggl|\frac{\Vert R_i Z_i-R_j Z_j \Vert
_2^2-\gamma
_{i,j}^2}{p} \biggr|.
\]
Using the fact that $a^2-b^2=(a-b)(a+b)$, and therefore, $|a^2-b^2|\leq
|a-b|(|a-b|+2|b|)$, we see that
\[
{\max_{i,j}}|a_{i,j}^2-b_{i,j}^2|\leq{\max_{i,j}}|a_{i,j}-b_{i,j}|
\Bigl({\max_{i,j}}|a_{i,j}-b_{i,j}|+{2\max}|b_{i,j}| \Bigr) .
\]
If we choose $a_{i,j}=\Vert R_i Z_i-R_j Z_j \Vert_2/\sqrt{p}$ and
$b_{i,j}=\gamma_{i,j}/\sqrt{p}$, we see that the previous equation becomes
\[
\Delta_1\leq\Delta_0 \biggl(\Delta_0+{2\max_{i\neq j}} \frac{\gamma
_{i,j}}{\sqrt{p}} \biggr) .
\]
Therefore, if we can show that $\Delta_0\max_{i\neq j}\gamma
_{i,j}/\sqrt
{p}$ goes to 0 in probability, we will have $\Delta_1\rightarrow0$ in
probability. Using the concentration result given in (\ref
{eq:ConditionalConcentration}), in connection with Proposition 1.9 in
\citet{ledoux2001} and a slight modification explained in
\citet{nekKernels}, we have
%
%
\begin{eqnarray}\label{eq:ControlGammaijSquare}
(R_i^2+R_j^2) \nu-\frac{\gamma_{i,j}^2}{p}&=&\operatorname{var}
_{\mathcal{YR}_n} \bigl(\Vert R_iZ_i-R_jZ_j \Vert_2/\sqrt{p}
\bigr)\nonumber\\[-8pt]\\[-8pt]
&\leq&
\frac{R^2_{\infty
}(p)}{p} \frac{32 C}{b(c_0)^{2/b}}\Gamma(2/b)=R^2_{\infty}(p)\frac
{\kappa_b}{p} .\nonumber
\end{eqnarray}
Using our assumption that $\nu$ remains bounded, we see that
\[
\frac{1}{R_{\infty}(p)}\max_{i\neq j}\frac{\gamma_{i,j}}{\sqrt
{p}}\qquad\mbox{remains bounded.}
\]
Therefore, for some $K$ independent of $p$,
\[
\max_{i\neq j}\frac{\gamma_{i,j}}{\sqrt{p}}\Delta_0\leq K
R_{\infty}(p)
r_1 ,
\]
with probability going to 1. Our assumptions also guarantee that
$R_{\infty}(p) r_1\rightarrow0$, so we conclude that, for a constant $K$
independent of $p$,
\begin{eqnarray}
{\max_{i\neq j}} \biggl|\frac{\Vert R_i Z_i-R_j Z_j \Vert_2^2-\gamma
_{i,j}^2}{p} \biggr|=\Delta_1\leq K r_1 R_{\infty}(p)
\rightarrow0\nonumber\\
\eqntext{\mbox{with probability going to 1}.}
\end{eqnarray}
Using (\ref{eq:ControlGammaijSquare}), we have the
deterministic inequality
\[
{\max_{i,j}} \biggl|(R_i^2+R_j^2) \nu-\frac{\gamma
_{i,j}^2}{p} \biggr|\leq R^2_{\infty}(p)\frac{\kappa_b}{p}\ll r_1^2
\ll
r_1 .
\]
So we can finally conclude that with high probability
\[
\Delta_2={\max_{i\neq j}}|\beta_{i,j}|={\max_{i\neq j}} \biggl|\frac
{\Vert R_i Z_i-R_j Z_j \Vert_2^2}{p}-\nu(R_i^2+R_j^2) \biggr|\leq K r_1
R_{\infty}(p)\rightarrow0 .
\]
Putting all these elements together, we see that when
\[
u_p=\frac{(\mathcal{M}_p^{1/2}\vee R_{\infty}(p))R_{\infty}(p)
(\log
n+(\log n)^{\varepsilon} )^{1/b}}{p^{1/2}} ,
\]
we can find a constant $K$ such that
\[
P\Bigl({\max_{i\neq j}}|2\alpha_{i,j}+\beta_{i,j}|>K u_p\Bigr)\rightarrow0 .
\]
In other words,
%
%
\begin{equation}\label{eq:ControlMaxErrorGeomApproxEntryWise}
P \Bigl({\max_{i\neq j}} \bigl|\Vert X_i-X_j \Vert_2^2- [\Vert
Y_i-Y_j \Vert_2^2+\nu(R_i^2+R_j^2) ] \bigr|>Ku_p
\Bigr)\rightarrow0 .
\end{equation}
This establishes (a strong form of) the first part of the theorem, that
is, (\ref{eq:InterpointDistanceEllipCase}).

- \textit{Second part of the theorem} [\textit{equation} (\ref
{eq:MainResConcEllipCase})].
To get to the second part, we recall that, assuming that $f$ is
$C_1$-Lipschitz on an interval containing $\{\Vert X_i-X_j \Vert
^2_2,\Vert Y_i-Y_j \Vert_2^2+\nu(R_i^2+R_j^2)\}$, we have
\begin{eqnarray*}
&&\bigl|f(\Vert X_i-X_j \Vert^2_2)-f\bigl(\Vert Y_i-Y_j \Vert_2^2+\nu
(R_i^2+R_j^2)\bigr) \bigr|\\
&&\qquad\leq C_1 \bigl|\Vert X_i-X_j \Vert^2_2-\bigl(\Vert Y_i-Y_j \Vert_2^2+\nu
(R_i^2+R_j^2)\bigr) \bigr|.
\end{eqnarray*}


Let us define, for $\eta>0$ given, the event
\[
E=\{\forall i\neq j, \Vert X_i-X_j \Vert_2^2 \in I_p(\eta) , \Vert
Y_i-Y_j \Vert_2^2 \in[\mathcal{W}_p,\mathcal{M}_p]\} ,
\]
and the random element
\[
\zeta_n=\sup_{f \in\mathcal{F}_{C_1,I_p(\eta)}}{\max_{i\neq
j}} \bigl|f(\Vert X_i-X_j \Vert^2_2)-f\bigl(\Vert Y_i-Y_j \Vert_2^2+\nu
(R_i^2+R_j^2)\bigr)
\bigr| .
\]

When $E$ is true, all the pairs $\{\Vert X_i-X_j \Vert^2_2,\Vert
Y_i-Y_j \Vert_2^2+\nu(R_i^2+R_j^2)\}$ are in $I_p(\eta)$: the
part concerning $\Vert Y_i-Y_j \Vert_2^2+\nu(R_i^2+R_j^2)$ is
obvious, and the one concerning $\Vert X_i-X_j \Vert^2_2$ comes from the
definition of $E$. So when $E$ is true, we also have
\[
\forall i\neq j \qquad \bigl|f(\Vert X_i-X_j \Vert^2_2)-f\bigl(\Vert
Y_i-Y_j \Vert_2^2+\nu(R_i^2+R_j^2)\bigr) \bigr|\leq C_1
|2\alpha
_{i,j}+\beta_{i,j} | .
\]
Let us now consider the random variable $\tau_n$ such that $\tau_n=C_1$
on $E$ and $\infty$ otherwise, so $\tau_n=C_11_{E}+\infty1_{E^c}$. Our
remark above shows that
\[
{\zeta_n\leq\tau_n \max_{i\neq j}} |2\alpha_{i,j}+\beta
_{i,j} |.
\]

Now, we see from our assumptions about $\{Y_i\}_{i=1}^n$,
(\ref{eq:ControlMaxErrorGeomApproxEntryWise}) and the fact that
$u_p\rightarrow0$, that for any $\eta>0$, $P(E)\rightarrow1$. So we have
\[
P(\tau_n\leq C_1)\rightarrow1 .
\]
Also, ${\max_{i\neq j}} |2\alpha_{i,j}+\beta_{i,j} |\leq K u_p$
with probability tending to 1, so we can conclude that
\[
P\Bigl({\tau_n \max_{i\neq j}} |2\alpha_{i,j}+\beta_{i,j}
|\leq C_1
K u_p\Bigr)\rightarrow1 .
\]
Hence, we also have
\[
P^* (\zeta_n \leq C_1K u_p )\rightarrow1 ,
\]
where this statement might have to be understood in terms of
outer\break
probabilities---hence the $P^*$ instead of $P$. [See \citet
{vandervaart}, page 258. In plain English, we have found a random
variable, ${\tau_n \max_{i\neq j} }|2\alpha_{i,j}+\beta
_{i,j}
|$, bounded by $C_1 K u_p$ with probability going to 1, which is larger
than the random element $\zeta_n$.]

In other respects, we have, for all $f \in\mathcal{F}_{C_1,I_p(\eta)}$,
\[
\Vert M_f-\widetilde{M}_f \Vert_F^2\leq\zeta_n^2 ,
\]
since
\begin{eqnarray*}
&&
{\max_{i,j}}|M_f(i,j)-\widetilde{M}_f(i,j)|\\
&&\qquad\leq{\frac{1}{n}\max
_{i\neq j}}
\bigl|f(\Vert X_i-X_j \Vert_2^2)-f\bigl(\Vert Y_i-Y_j \Vert_2^2+\nu
(R_i^2+R_j^2)\bigr)\bigr|\leq\frac{\zeta_n}{n} .
\end{eqnarray*}
Therefore,
%
%
\begin{equation}\label{eq:keyEqControlSupErrorEllipConcCaseProof}
{\sup_{f \in\mathcal{F}_{C_1,I_p(\eta)}}}\Vert M_f-\widetilde{M}_f
\Vert_F\leq
\zeta_n
\rightarrow0 \qquad\mbox{in probability},
\end{equation}
where once again this statement may have to be understood in terms of
outer probabilities. The result stated in (\ref
{eq:MainResConcEllipCase}) is proved.
\end{pf*}

We mentioned before the proof the possibility that we might let $C_1$
vary with $n$ and $p$ and still get a good approximation result. This
can be done by looking at (\ref
{eq:keyEqControlSupErrorEllipConcCaseProof}) above: $\zeta_n$ is less
than $KC_1 u_p$ with high probability, so when $u_p C_1(n)\rightarrow0$
the main approximation result of Theorem
\ref{thm:InfoPlusNoiseConcCase} holds, for a $C_1$ and therefore a class of
functions, that vary with $n$ (and $p$).

\subsubsection{The case of inner-product kernel random matrices}

We now turn our attention to kernel matrices of the form
$M(i,j)=f(X_i^{\prime}X_j)/n$ which are also of interest in practice. In
that setting, we are able to obtain results similar in flavor to
Theorem \ref{thm:InfoPlusNoiseConcCase}, with slight modifications on
the assumptions we make about $f$.
\begin{theorem}[(Scalar product kernels)]\label
{thm:InfoPlusNoiseConcCaseDotProductsKernels}
Suppose we observe data $X_1,\ldots,\break X_n$ in $\mathbb{R}^p$, with
\[
X_i=Y_i+R_i \frac{Z_i}{\sqrt{p}} .
\]
We place ourselves in the high-dimensional setting where $n$ and $p$
tend to infinity.
We assume that $\{Y_i\}_{i=1}^n\sim P_n$.

$\{Z_i\}_{i=1}^n$ are i.i.d. with $\mathbf{E}(Z_i)=0$, and we also assume
that $\{Y_i\}_{i=1}^n $ and $\{Z_i\}_{i=1}^n$ are independent.

$\{R_i\}_{i=1}^n$ are assumed to be independent of $\{Z_i\}_{i=1}^n$.
We also assume that we can find a deterministic sequence $R_{\infty
}(p)$ such that $\forall i, |R_i|\leq R_{\infty}(p)$ and $R_{\infty
}(p)\geq1$.

We assume that the distribution of $Z_i$ is such for any 1-Lipschitz
function $F$ (with respect to Euclidean norm), if $\mu_F=\mathbf{E}(F(Z_i))$,
\[
P\bigl(|F(Z_i)-\mu_F|>r\bigr)\leq C \exp(-c_0 r^b)\triangleq h(r) ,
\]
where for simplicity we assume that $c_0$, $C$ and $b$ are independent
of $p$. We call $\nu=\mathbf{E}(\Vert Z_i \Vert_2^2)/p$ and assume that
$\nu$ stays bounded as $p\rightarrow\infty$.

We call $\mathcal{M}={\max_{i,j}} |Y_i^{\prime}Y_j |$, and
$\mathcal{M}_p$ a real such that $P(\mathcal{M}\leq\mathcal
{M}_p)\rightarrow1$. We
assume that there exists $\varepsilon>0$ such that
\[
\max(\mathcal{M}_p^{1/2},R_{\infty}(p)) \frac{R_{\infty}(p)(\log n
+(\log
n)^{\varepsilon})^{1/b}}{\sqrt{p}}\rightarrow0.
\]

We then have
%
%
\begin{equation}\label{eq:dotProductsEllipCase}
{\max_{i,j}} |X_i^{\prime}X_j-(Y_i^{\prime}Y_j+\delta
_{i,j}\nu
R_i^2) |\rightarrow0 \qquad\mbox{in probability}.
\end{equation}

We call $J_p(\eta)=[-\mathcal{M}_p-\eta-R_{\infty
}^2(p)\nu,\mathcal{M}_p+\eta+R_{\infty}^2(p)\nu
]$ and
\[
\mathcal{F}_{C_1,J_p(\eta)}=\Bigl\{f \mbox{ such that }
{\sup_{x,y\in J_p(\eta)}} |f(x)-f(y)|\leq C_1|x-y|\Bigr\}.
\]

We then consider the random matrices $M_f$ with $(i,j)$ entry
\[
M_f(i,j)=\frac{1}{n}f (X_i^{\prime}X_j ) \qquad \mbox{for } f
\in
\mathcal{F}_{C_1,J_p(\eta)} .
\]

Let us call $\widetilde{M}$ the matrix with $(i,j)$th entry
\[
\widetilde{M}_f(i,j)=
\cases{
\dfrac{1}{n}f (Y_i^{\prime}Y_j ), &\quad if $i\neq j$,\vspace*{2pt}\cr
\dfrac{1}{n} f(\Vert Y_i \Vert_2^2+\nu R_i^2), &\quad if $i=j$.}
\]
We have, for any $C_1>0$ and $\eta>0$,
\[
{\lim_{n,p\rightarrow\infty} \sup_{f\in\mathcal{F}_{C_1,J_p(\eta)}}}
\Vert M_f-\widetilde{M}_f \Vert_F= 0 \qquad\mbox{in probability.}
\]
\end{theorem}

We note that under our assumptions, we also have $|f(\Vert Y_i \Vert
_2^2+\nu R_i^2)-f(\Vert Y_i \Vert_2^2)|\leq\nu
C_1R_{\infty}^2(p)$, with high probability, and uniformly in $f$ in
$\mathcal{F}_{C_1,J_p(\eta)}$.
Therefore, when\vspace*{1pt} $R_{\infty}^4(p)/n\rightarrow0$, the result is also valid
if we replace the diagonal of $\widetilde{M}_f$ by $\{f(\Vert Y_i
\Vert_2^2)\}
_{i=1}^n/n$---in which case the new approximating matrix is the kernel
matrix computed from the signal part of the data. Furthermore, the same
argument shows that we get a valid operator norm approximation of $M$
by this ``pure signal'' matrix as soon as $R_{\infty}^2(p)/n$ tends to 0.

The same measurability issues as in the previous theorems might arise
here and the statement should be understood as before: we can find a
random variable going to 0 in probability that is larger than the
random element ${\sup_{f\in\mathcal{F}_{C_1,J_p(\eta)}}} \Vert
M_f-\widetilde{M}_f \Vert_F$.

Finally, let us note that once again the theorem is stated for a fixed
$C_1$ [and hence for an essentially fixed (with $n$) class of
functions, though some changes in this class might come from varying
$J_p(\eta)$], but the proof allows us to deal with a varying $C_1(n)$.
The adjustments are very similar to the ones we discussed after the
proof of Theorem \ref{thm:InfoPlusNoiseConcCase} and we leave them to
the interested reader.
\begin{pf*}{Proof of Theorem \ref{thm:InfoPlusNoiseConcCaseDotProductsKernels}}
The proof is quite similar to that of Theorem
\ref{thm:InfoPlusNoiseConcCase}, so we mostly outline the differences and
use the same notation as before. We now have to focus on
\[
X_i^{\prime}X_j=Y_i^{\prime}Y_j +R_i\frac{Z_i^{\prime}Y_j}{\sqrt
{p}}+R_j\frac
{Z_j^{\prime}Y_i}{\sqrt{p}}+R_iR_j\frac{Z_i^{\prime}Z_j}{p} .
\]
The analysis of $R_i\frac{Z_i^{\prime}Y_j}{\sqrt{p}}$ is entirely similar
to our analysis of $\alpha_{i,j}$ in the proof of Theorem
\ref{thm:InfoPlusNoiseConcCase}. The key remark now is that as function of
$Z_i$, when $\mathcal{YR}_n\in\mathcal{LR}$, it is, with the new definition
of $\mathcal{M}_p$, $R_{\infty}(p)\sqrt{\mathcal{M}_p/p}$-Lipschitz with
respect to Euclidean norm. So we immediately have, with the new
definition of $\mathcal{M}_p$: if $r_0=R_{\infty}(p)(\mathcal
{M}_p/p)^{1/2}(\log n+(\log n)^{\varepsilon})^{1/b}(2/c_0)^{1/b}$, and
$\mathcal{YR}_n \in\mathcal{LR}$, for some $K>0$ which does not
depend on $\mathcal{YR}_n$,
\[
P_{\mathcal{YR}_n} \biggl({\max_{i,j}} \biggl|R_i\frac{Z_i^{\prime}
Y_j}{\sqrt
{p}} \biggr|>r_0 \biggr)\leq K \exp(-2(\log n)^{\varepsilon}) .
\]
Now, since $P(\mathcal{YR}_n\notin\mathcal{LR})\rightarrow0$, we
conclude as
before that
\[
P \biggl({\max_{i,j}} \biggl|R_i\frac{Z_i^{\prime}Y_j}{\sqrt{p}}
\biggr|>r_0 \biggr)\rightarrow0 .
\]

On the other hand, using the fact that $4 R_iR_j Z_i^{\prime}Z_j=\Vert
R_iZ_i+R_jZ_j \Vert_2^2-\Vert R_iZ_i-R_jZ_j \Vert_2^2$, and analyzing the
concentration properties of $\Vert R_iZ_i+R_jZ_j \Vert_2^2$ in the
same way
as we did those of $\Vert R_iZ_i-R_jZ_j \Vert_2^2$, we conclude that if
$u_p=R_{\infty}^2(p)(2/c_0)^{1/b}(\log n+(\log n)^{\varepsilon
})^{1/b}p^{-1/2}$, we can find a constant $K$ such that
\[
P \biggl({\max_{i\neq j}} \biggl|\frac{\Vert R_iZ_i-R_jZ_j \Vert
_2^2}{p}-\nu(R_i^2+R_j^2) \biggr|>K u_p
\biggr)\rightarrow0
\]
and
\[
P \biggl({\max_{i\neq j}} \biggl|\frac{\Vert R_iZ_i+R_jZ_j \Vert
_2^2}{p}-\nu(R_i^2+R_j^2) \biggr|>K u_p
\biggr)\rightarrow0 .
\]
Similar arguments, relying on the fact that \mbox{$\Vert\cdot\Vert_2$} is
obviously 1-Lipschitz with respect to Euclidean norm, also lead to the
fact that
\[
P \biggl({\max_{i}} \biggl|\frac{R_i^2\Vert Z_i \Vert_2^2}{p}-\nu
R_i^2 \biggr|>Ku_p \biggr)\rightarrow0 .
\]
Therefore, we can find $K$, greater than 1 without loss of generality,
such that
\[
P \biggl({\max_{i,j}} \biggl|R_iR_j\frac{Z_i^{\prime}Z_j}{p}-\delta
_{i,j}\nu R_i^2 \biggr|>K u_p \biggr)\rightarrow0 .
\]
We can therefore conclude that
\[
P \Bigl({\max_{i,j}} |X_i^{\prime}X_j-(Y_i^{\prime}Y_j+\delta
_{i,j}\nu R_i^2) |>Ku_p+2r_0 \Bigr)\rightarrow0 .
\]

If $R_{\infty}(p)\max((\mathcal{M}_p)^{1/2},R_{\infty}(p)) (\log
n+(\log
n)^{\varepsilon})^{1/b}/\sqrt{p}\rightarrow0$, then both $r_0$ and
$u_p$ tend to
0. Therefore, under our assumptions,
\[
{\max_{i,j}} |X_i^{\prime}X_j-(Y_i^{\prime}Y_j+\delta
_{i,j}\nu
R_i^2) |\rightarrow0 \qquad\mbox{in probability.}
\]
So we have shown the first assertion of the theorem.

The final step of the proof is now clear: we have, for all $(i,j)$,
\[
|f(X_i^{\prime}X_j)-f(Y_i^{\prime}Y_j+\delta_{i,j}\nu
R_i^2) |\leq C_1 |X_i^{\prime}X_j-(Y_i^{\prime}Y_j+\delta
_{i,j}\nu R_i^2) | ,
\]
when for all $(i,j)$, $X_i^{\prime}X_j$ and $(Y_i^{\prime}Y_j+\delta
_{i,j}\nu R_i^2)$ are in $J_p(\eta)$. This event happens with
probability going to 1 under our assumptions. So following the same
approach as before and dealing with measurability in the same way, we
have, with probability going to 1,
\begin{eqnarray*}
&&
{\sup_{f\in\mathcal{F}_{C_1,J_p(\eta)}}\max_{i,j}}
|f(X_i^{\prime}
X_j)-f(Y_i^{\prime}Y_j+\delta_{i,j}\nu R_i^2) |
\\
&&\qquad\leq {C_1
\max_{i,j}} |X_i^{\prime}X_j-(Y_i^{\prime}Y_j+\delta_{i,j}\nu
R_i^2) | .
\end{eqnarray*}
So we conclude that
\[
{\sup_{f\in\mathcal{F}_{C_1,J_p(\eta)}}\max_{i,j}}
|f(X_i^{\prime}
X_j)-f(Y_i^{\prime}Y_j+\delta_{i,j}\nu R_i^2) | \rightarrow0
\qquad\mbox{in probability.}
\]
From this statement, we get in the same manner as before,
\[
\sup_{f\in\mathcal{F}_{C_1,J_p(\eta)}}\Vert M_f-\widetilde{M}_f
\Vert_F
\rightarrow0\qquad
\mbox{in probability.}
\]
\upqed\end{pf*}

As before, the equations above show that if $C_1(n)(u_p+r_0)\rightarrow
0$, the same approximation result holds, now with a varying $C_1(n)$.

\subsection{Practical consequences of the results: Case of spherical
noise}\label{subsec:practicalConsequences}

Our aim in giving approximation results is naturally to use existing
knowledge concerning the approximating matrix to reach conclusions
concerning the information${}+{}$noise kernel matrices that are of interest
here. In particular, we have in mind situations where the ``signal''
part of the data, that is, what we called $\{Y_i\}_{i=1}^n$ in the
theorems, and $f$ [or $f(\cdot+2\nu)$, with $\nu$
being as defined in Theorems \ref{thm:InfoPlusNoiseGaussianCase} or
\ref{thm:InfoPlusNoiseConcCase}] are such that the assumptions of Theorems
3.1 or 5.1 in \citet{KoltchinGine00} are satisfied, in which case
we can approximate the eigenvalues of $\widetilde{M}$ by those of the
corresponding operator in $L^2(dP)$. In this setting the matrix
$\widetilde{M}
$, which is normalized so its entries are of order $1/n$ has a
nondegenerate limit, which is why we considered for our kernel matrices
the normalization $f(\Vert X_i-X_j \Vert_2^2)/n$. [This normalization by
$1/n$ makes our proofs considerably simpler than the ones given in
\citet{nekKernels}.]

Another potentially interesting application is the case where the
signal part of the data is sampled i.i.d. from a manifold with bounded
Euclidean diameter, in which case our results are clearly applicable.

\subsubsection{Spectral properties of information${}+{}$noise kernel random
matrices from pure signal kernel random matrices}

The practical interest of the theorems we obtained above lie in the
fact that the Frobenius norm is larger than the operator norm, and
therefore all of our results also hold in operator norm. Now we recall
the discussion in El Karoui [(\citeyear{nekSparseMatrices}), Section 3.3], where we
explained that consistency in operator norm implies consistency of
eigenvalues and consistency of eigenspaces corresponding to separated
eigenvalues [as consequences of Weyl's inequality and the Davis--Kahane
$\sin(\theta)$ theorem---see \citet{bhatia97} and \citet
{stewart90}].

Theorems \ref{thm:InfoPlusNoiseGaussianCase}, \ref
{thm:InfoPlusNoiseConcCase}, \ref
{thm:InfoPlusNoiseConcCaseDotProductsKernels} therefore imply that
under the assumptions stated there, the spectral properties of the
matrix $M$ can be deduced from those of the matrix $\widetilde{M}$. In
particular, for techniques such as kernel PCA, we expect, when it is a
reasonable idea to use that technique, that $M$ will have some
separated eigenvalues, that is, a few will be large and there will be a
gap in the spectrum. In that setting, it is enough to understand
$\widetilde{M}$, which corresponds, if $\forall i, R_i=1$, to a pure signal
matrix, with a possibly slightly different kernel, to have a
theoretical understanding of the properties of the technique.

For instance, if $\forall i, R_i=1$, if the assumptions underlying the
first-order results of \citet{KoltchinGine00} are satisfied for
$\widetilde{M}$, the (first-order) spectral properties of $M$ are the
same as
those of $\widetilde{M}$, and hence of the corresponding operator in $L^2(dP)$.

\subsubsection{On the Gaussian kernel}

Our analysis reveals a very interesting feature of the Gaussian kernel,
that is, the case where $M(i,j)=\exp(-s\Vert X_i-X_j \Vert_2^2)/n$,
for some
$s>0$: when Theorem \ref{thm:InfoPlusNoiseGaussianCase} or Corollary
\ref{coro:SphericalNoiseKernels} (i.e., Theorem \ref
{thm:InfoPlusNoiseConcCase} with $\forall i, R_i=1$) apply, the
eigenspaces corresponding to separated eigenvalues of the signal${}+{}$noise
kernel matrix converge to those of the pure signal matrix.

This is simply due to the fact that in that setting, if $\mathcal{S}$ is
the matrix such that
\[
\mathcal{S}(i,j)=\exp(-2\nu s) \frac{1}{n} \exp(-s\Vert Y_i-Y_j
\Vert_2^2) ,
\]
a rescaled version of the ``pure signal'' matrix $\mathcal{M}$ with
$(i,j)$th entry $\frac{1}{n} \exp(-s\Vert Y_i-Y_j \Vert_2^2)$, we have
\[
\interleave\mathcal{S}-\widetilde{M}\interleave_2\rightarrow0 .
\]
This latter statement is a simple consequence of the fact that
$\mathcal{S}-\widetilde{M}$ is a diagonal matrix with entries $(\exp
(-2\nu
s)-1)/n$ on the diagonal, and therefore its operator norm goes to 0. On
the other hand, $\mathcal{S}$ clearly has the same eigenvectors as the
pure signal matrix $\mathcal{M}$. Hence, because the eigenspaces of
$\widetilde{M}$ are consistent for the eigenspaces of $\mathcal{S}$
corresponding to separated eigenvalues, they are also consistent for
those of $\mathcal{M}$. (We note that our results are actually stronger
and allow us to deal with a collection of matrices with varying $s$ and
not a single $s$, as we just discussed. This is because we can deal
with approximations over a collection of functions in all our theorems.)

Because of the practical importance of eigenspaces in techniques such
as kernel PCA, these remarks can be seen as giving a theoretical
justification for the use of the Gaussian kernel over other kernels in
the situations where we think we might be in an information${}+{}$noise
setting, and the noise is spherical.

On the other hand, $\mathcal{S}$ underestimates the large eigenvalues of
$\mathcal{M}$ because $\mathcal{S}=\exp(-2\nu s)\mathcal
{M}$, and
obviously $\exp(-2\nu s)<1$. Using Weyl's inequality [see
\citet{bhatia97}], we have, if we denote by $\lambda_i(M)$ is the
$i$th eigenvalue of the symmetric matrix $M$,
\[
\forall i , 1\leq i \leq n, \qquad |\lambda_i(\widetilde
{M})-\lambda
_i(\mathcal{S}) |\leq\interleave\widetilde{M}-\mathcal{S}
\interleave_2 .
\]
Since the right-hand side goes to 0 asymptotically, the eigenvalues of
$\mathcal{M}$ (the ``pure signal'' matrix) that stay asymptotically\vspace*{1pt} bounded
away from 0 are underestimated by the corresponding eigenvalues of
$\widetilde{M}$.

When the noise is elliptical, that is, $R_i$'s are not all equal to 1,
the ``new'' matrix $\mathcal{S}$ we have to deal with has entries
\[
\mathcal{S}(i,j)=\exp(-sR_i^2)\exp(-sR_j^2)\frac{1}{n} \exp
(-s\Vert Y_i-Y_j \Vert_2^2) ,
\]
so it can be written in matrix form
\[
\mathcal{S}=D\mathcal{M} D ,
\]
where $D$ is a diagonal\vspace*{1pt} matrix with $D(i,i)=\exp(-sR_i^2)$. By the same
arguments as above, $\interleave S-\widetilde{M}\interleave_2
\rightarrow0$ in
probability, but
now $\mathcal{S}$ does not have the same eigenvectors as the pure signal
matrix $\mathcal{M}$. So in this elliptical setting if we were to do
kernel analysis on $M$, we would not be recovering the eigenspaces of
the pure signal matrix $\mathcal{M}$.

\subsubsection{Variants of kernel matrices: Laplacian matrices and the
issue of centering}

In various parts of statistics and machine learning, it has been argued
that Laplacian matrices should be used instead of kernel matrices. See,
for instance, the very interesting \citet
{BelkinNiyogiCvLapEigenmaps08}, where various spectral properties of
Laplacian matrices have been studied, under a ``pure'' signal assumption
in our terminology. For instance, it is assumed that the data is
sampled from a fixed-dimensional manifold. In light of the theoretical
and practical success of these methods, it is natural to ask what
happens in the information${}+{}$noise case.

There are several definitions of Laplacian matrices. A popular one
[see, e.g., the work of \citet{BelkinNiyogiCvLapEigenmaps08}, among other publications],
is derived from kernel
matrices: given $M$ a kernel matrix, the Laplacian matrix is defined as
\[
L(i,j)=
\cases{
-M(i,j),&\quad if $i\neq j$,\vspace*{2pt}\cr
\displaystyle\sum_{i\neq j} M(i,j), &\quad otherwise.}
\]

When our Theorems \ref{thm:InfoPlusNoiseConcCase} or \ref
{thm:InfoPlusNoiseConcCaseDotProductsKernels} apply, we have seen that,
for relevant classes of functions $\mathcal{F}$, ${\sup_{f\in\mathcal{F}}
n\max_{i\neq j}}|M_f(i,j)-\widetilde{M}_f(i,j)|\rightarrow0$ in probability.

Let us now focus on the case of a single function $f$. If we call
$\widetilde{L}$ the Laplacian matrix corresponding to $\widetilde
{M}$, we have
\begin{eqnarray*}
{n\max_{i\neq j}}|L(i,j)-\widetilde{L}(i,j)|&\rightarrow&0 \qquad\mbox{in
probability},\\
{\max_{i}}|L(i,i)-\widetilde{L}(i,i)|&\rightarrow&0 \qquad\mbox{in
probability}.
\end{eqnarray*}

We conclude that $\interleave L-\widetilde{L}\interleave_2
\rightarrow0$ in probability;
we can therefore deduce that the spectral properties of the Laplacian
matrix $L$ from those of $\widetilde{L}$, which, when $\forall i,
R_i=1$, is a ``pure signal'' matrix, where we have slightly adjusted the
kernel. Here again, the Gaussian kernel plays a special role, since
when we use a Gaussian kernel, $\widetilde{L}$ is a scaled version of
the Laplacian matrix computed from the signal part of the data.

Finally, other versions of the Laplacian are also used in practice. In
particular, a ``normalized'' version is sometimes advocated, and
computed as $N_L=D_L^{-1/2}L D_L^{-1/2}$, if $D$ is the diagonal of the
matrix $L$ defined above. We have just seen that $\interleave
D_L-D_{\widetilde{L}}\interleave_2\rightarrow0$ in probability
and $\interleave L-\widetilde{L}\interleave_2\rightarrow0$ in
probability. Therefore, if the entries
of $D_{\widetilde{L}}$ are bounded away from 0 with probability going
to~1, we conclude that $\interleave D^{-1}_{\widetilde{L}}
\interleave_2$ stays bounded
with high probability and
\[
\interleave N_L-N_{\widetilde{L}}\interleave_2\rightarrow0 \qquad\mbox
{in probability.}
\]
So once again, understanding the spectral properties of $N_L$
essentially boils down to understanding those of $N_{\widetilde{L}}$,
which is, in the spherical setting where $\forall i, R_i=1$, a ``pure
signal'' matrix. In the case of the Gaussian kernel, $N_{\widetilde{L}}$
is equal to the normalized Laplacian matrix computed from the ``pure
signal'' data $\{Y_i\}_{i=1}^n$.

\paragraph*{The question of centering} In practice, it is often the case
that one works with centered versions of kernel matrices: either the
row sums, the column sums or both are made to be equal zero. These
centering operations amount to multiplying (resp., on the right,
left or both) our original kernel matrix by the matrix $H=\operatorname
{Id}_n-\mathbf{11}^{\prime}/n$, where $\mathbf{1}$ is the $n$-dimensional vector whose
entries are all equal to 1. This matrix has operator norm 1, so when
$\widetilde{M}$ is such that $\interleave M-\widetilde{M}
\interleave_2\rightarrow
0$, the same is true
for $H^aMH^b$ and $H^a\widetilde{M}H^b$, where $a$ and $b$ are either
0 or
1. This shows that our approximations are therefore also informative
when working with centered kernel matrices.

\section{Conclusions}\label{sec3}
Our results aim to bridge the gap in the existing literature between
the study of kernel random matrices in the presence of pure
low-dimensional signal data [see, e.g., \citet{KoltchinGine00}]
and the case of truly high-dimensional data [see \citet
{nekKernels}]. Our study of information${}+{}$noise kernel random matrices
shows that, to first order, kernel random matrices are somewhat
``spectrally robust'' to the corruption of signal by additive
high dimensional and spherical noise (whose norm is controlled). In
particular, they tend to behave much more like a kernel matrix computed
from a low-dimensional signal than one computed from high-dimensional
data.

Some noteworthy results include the fact that dot-product kernel random
matrices are, under reasonable assumptions on the kernel and the
``signal distribution'' spectrally robust for both eigenvalues and
eigenvectors. The Gaussian kernel also yields spectrally robust
matrices at the level of eigenvectors, when the noise is spherical.
However, it will underestimate separated eigenvalues of the Gaussian
kernel matrix corresponding to the signal part of the data.

On the other hand, Euclidean distance kernel random matrices are not,
in general, robust to the presence of additive noise. As our results
show, under reasonably minimal assumptions on both the noise, the
kernel and the signal distribution, a Euclidean distance kernel random
matrix computed from additively corrupted data behaves like another
Euclidean distance kernel matrix computed from another kernel: in the
case of spherical noise, it is a shifted version of $f$, the shift
being twice the norm of the noise. For spherical noise, this is bound
to create (except for the Gaussian kernel) potentially serious
inconsistencies in both estimators of eigenvalues and eigenvectors,
because the eigenproperties of the kernel matrix corresponding to the
function $f_{\nu}(\cdot)=f(\cdot+2\nu)$ are in
general different from that of the kernel matrix corresponding to the
function $f$. The same remarks apply to the case of elliptical noise,
where the change of kernel is not deterministic and even more
complicated to describe and interpret.

Our study also highlights the importance of the implicit geometric
assumptions that are made about the noise. In particular, the results
are qualitatively different if the noise is spherical (e.g.,
multivariate Gaussian) or elliptical (e.g., multivariate~$t$).
Interpretation is more complicated in the elliptical case and a number
of nice properties (e.g., robustness or consistency) which hold for
spherical noise do not hold for elliptical noise.

We note that our study suggests that simple practical (and entrywise)
corrections could be used to go from the ``signal${}+{}$noise'' situation to
an approximation of the ``pure signal'' situation. Those would naturally
depend on the noise geometry and what information practitioners have
about it.

Our results can therefore be seen as highlighting (from a theoretical
point of view) the strength and limitations of techniques which rely on
kernel random matrices as a primary element in a data analysis. We hope
they shed light on an interesting issue and will help refine our
understanding of the behavior of kernel techniques and related
methodologies for high-dimensional input data.


\section*{Acknowledgments}
The author would like to thank Peter Bickel for suggesting that he
consider the problem studied here and in general for many enlightening
discussions about topics in high-dimensional statistics. He would also
like to thank an anonymous referee and Associate Editors for raising
interesting questions which contributed to a strengthening of the
results presented in the paper.

\printaddresses

\end{document}